\documentclass[11pt,english]{article}
\usepackage{lmodern}

\usepackage[T1]{fontenc}
\usepackage[a4paper]{geometry}
\usepackage{color}
\definecolor{document_fontcolor}{rgb}{0, 0, 0}
\color{document_fontcolor}
\usepackage{textcomp}
\usepackage{mathrsfs}
\usepackage{amsmath}
\usepackage{amsthm}
\usepackage{amssymb}

\makeatletter
\newcommand{\lyxaddress}[1]{
\par {\raggedright #1
\vspace{1.4em}
\noindent\par}
}
\theoremstyle{plain}
\ifx\thechapter\undefined
\newtheorem{thm}{\protect\theoremname}
\else
\newtheorem{thm}{\protect\theoremname}[chapter]
\fi
\theoremstyle{definition}
\newtheorem{defn}[thm]{\protect\definitionname}
\theoremstyle{plain}
\newtheorem{prop}[thm]{\protect\propositionname}
\ifx\proof\undefined
\newenvironment{proof}[1][\protect\proofname]{\par
\normalfont\topsep6\p@\@plus6\p@\relax
\trivlist
\itemindent\parindent
\item[\hskip\labelsep\scshape #1]\ignorespaces
}{%
\endtrivlist\@endpefalse
}
\providecommand{\proofname}{Proof}
\fi
\theoremstyle{remark}
\newtheorem{rem}[thm]{\protect\remarkname}
\theoremstyle{plain}

\makeatother

\usepackage{babel}
\providecommand{\corollaryname}{Corollary}
\providecommand{\definitionname}{Definition}
\providecommand{\propositionname}{Proposition}
\providecommand{\remarkname}{Remark}
\providecommand{\theoremname}{Theorem}
\date{}

\begin{document}

\title{Polarized vectorial Poisson structures}

\author{Azzouz AWANE. Ismail BENALI. Souhaila EL AMINE}
\maketitle

\lyxaddress{LAMS. Ben M'sik's Faculty of Sciences. B.P.7955. Bd Driss Harti. Casablanca. Hassan II University of Casablanca.\\ 
 azzouz.awane@univh2c.ma; ismail.benali-etu@etu.univh2c.ma; souhaila.elamine-etu@etu.univh2c.ma
}
\begin{abstract}
We study various properties of polarized vectorial
Poisson  structures subordinate to polarized $k$-symplectic manifolds,
and also, we study the notion of polarized vectorial Poisson manifold.
Some properties and examples are given. 

M.S.C2010: 53D05, 37K05, 53D17.

 Hamiltonian. polarized vectorial Poisson manifold.
\end{abstract}

\maketitle

\section{{\large{}Introduction}}

\noindent Mathematical and physics considerations have led to introduce
the polarized $k-$symplectic structures (\cite{key-AW}, \cite{key-AW-GZ}
and \cite{key-NB}). The Poisson aspect of polarized $k$-symplectic
manifolds allows us to introduce and study the notion of polarized
vectorial Poisson structures. 

\noindent Recall that a Polarized $k$-symplectic structure on an
$n(k+1)$-dimensional foliated manifold $M$ is a pair $\left(\theta,\mathfrak{F}\right)$
in which $\mathfrak{F}$ is an $n-$codimensional foliation and $\theta$
is a closed and nondegenerate $\mathbb{R}^{k}$-valued differential
$2-$form vanishing on vector fields tangent to the leaves of $\mathfrak{F}$.

\noindent The polarized $k$-symplectic Darboux's theorem show that
around each point $x_{0}$ of $M$ there is a local coordinate system
$(x^{pi},y^{i})_{1\leq p\leq k,1\leq i\leq n}$ such that 
\[
\theta=\sum_{p=1}^{k}\left(\sum_{i=1}^{n}dx^{pi}\wedge dy^{i}\right)\otimes v_{p}
\]
and, $\mathfrak{F}$ is defined by the equations : $dy^{1}=0,\ldots,dy^{n}=0.$
Where $\left(v_{p}\right)_{1\leq p\leq k}$ is the canonical basis
of $\mathtt{\mathbb{R}}^{\mathit{k}}.$

A polarized Hamiltonian vector field is a foliate vector field $X$
such that $i\left(X\right)\theta$ is exact. An associated polarized
Hamiltonian to $X$ is an $\mathbb{R}^{k}$-valued function $H\in\mathfrak{\mathscr{C^{\infty}}}\left(M,\mathbb{R^{\mathit{k}}}\right)$
such that $i\left(X\right)\theta=-dH$. Locally the polarized Hamiltonians
 have the following form 
\[
H=\sum_{p=1}^{k}\left(\sum_{j}a_{j}(y^{1},...,y^{n})x^{pj}+b^{p}(y^{1},...,y^{n})\right)\otimes v_{p}
\]

where $a_{j}$ and $b_{p}$ are basic functions. 

The set of all polarized Hamiltonians is a proper vector subspace
of $\mathscr{C}^{\infty}(M,\mathbb{R}^{k})$, that we denote by $\mathfrak{H}(M,\mathfrak{F})$.
This subspace admits a natural Lie algebra law $\left\{ ,\right\} $,
called polarized vectorial Poisson structure subordinate to the polarized
$k-$symplectic structure.

In this work, we study various aspects of polarized vectorial Poisson
manifolds and we give some properties and examples of polarized Hamiltonians.
This leads us to introduce the notion of polarized vectorial Poisson
structure on a foliated manifold. We give in this paper, a natural
polarized $k-$symplectic structure $\left(\theta,\mathfrak{F}\right)$
on the space $\hom\left(\mathcal{G},\mathbb{R}^{k+1}\right)$, for
a given real Lie algebra $\mathcal{G}$; and also, the associated
linear polarized Poisson structure on $\mathfrak{H}\left(\hom\left(\mathcal{G},\mathbb{R}^{k+1}\right),\mathfrak{F}\right)$
having for support the space $\mathfrak{H}\left(\hom\left(\mathcal{G},\mathbb{R}^{k+1}\right),\mathfrak{F}\right)$
depending on the Lie algebra law of $\mathcal{G}$.

\section{Polarized $k$-symplectic manifolds}

Let $M$ be an $n(k+1)$-dimensional smooth manifold endowed with
an $n$-codimensional foliation $\mathfrak{F}.$ Let $\theta=\theta^{p}\otimes v_{p}\in\mathscr{A}_{2}(M)\otimes\mathbb{R}^{k}$
be an $\mathbb{R}^{k}$-valued differential $2-$form. We denote by
$E$ the sub-bundle of $TM$ defined by the tangent vectors of the
leaves of the foliation $\mathfrak{F}$. And also, we denote by $\Gamma(E)$,
the set of all cross-sections of the $M$-bundle $E\longrightarrow M$,
and by $\mathscr{A}_{p}(M)$ the set of all differential $p$-forms
on $M$. $\left(v_{p}\right)_{1\leq p\leq k}$ being the canonical
basis of the real vector space $\mathbb{R}^{k}$.

We recall that (\cite{key-AW}, \cite{key-AW-GZ}), $\left(\theta,E\right)$
is a polarized $k$-symplectic structure on $M$ if: (i) $\theta$
closed i.e. $d\theta=0$ ; (ii) $\theta$ nondegenerate, i.e., for
all $X\in\mathfrak{X}\left(M\right)$, $i\left(X\right)\theta=0\Longrightarrow X=0$
and (iii) $\theta(X,Y)=0$ for all $X,Y\in\Gamma\left(E\right).$

We recall also the following theorem (\cite{key-AW}, \cite{key-AW-GZ}),
which gives the local model of a polarized $k$-symplectic structure
in the Darboux's sense.
\begin{thm}
If $\left(\theta,E\right)$ is a polarized $k$-symplectic structure
on $M$, then for every point $x_{0}$ of $M$, there exist an open
neighborhood $U$ of $M$ containing $x_{0}$ equipped with a local
coordinate system $(x^{pi},y^{i})_{1\leq p\leq k,1\leq i\leq n}$
called an adapted coordinate system, such that the $\mathbb{R}^{k}$-valued
differential $2-$form $\theta$ is represented on $U$ by 
\[
\theta_{\mid U}=\sum_{p=1}^{k}\theta^{p}\otimes v_{p}=\sum_{p=1}^{k}\left(\sum_{i=1}^{n}dx^{pi}\wedge dy^{i}\right)\otimes v_{p}
\]
 and $\mathfrak{F}_{\mid U}$ is defined by the equations $dy^{1}=...=dy^{n}=0.$ 

The theorem's expressions imply the following local transition formulas
of the canonical coordinates 
\[
\overline{y}^{i}=\overline{y}^{i}\left(y^{1},\ldots,y^{n}\right),\,\,\,\,\,\,\,\,\overline{x}^{pi}=\sum_{j=1}^{n}x^{pj}\frac{\partial y^{j}}{\partial\overline{y}^{i}}+\varphi^{pi}\left(y^{1},\ldots,y^{n}\right)
\]
\end{thm}

Indeed, these expressions are affine with respect to $x^{pj}$.

Recall that, \cite{key-MLN2}, a real function $f\in\mathscr{C}^{\infty}\left(M\right)$
is called basic, if for any vector field $Y$ tangent to $\mathfrak{F}$,
the function $Y(f)$ is identically zero. We denote by $\mathscr{A}_{b}^{0}(M,\mathfrak{F})$
the subring of $\mathscr{C}^{\infty}\left(M\right)$ of basic functions. 

Let $f\in\mathscr{C}^{\infty}\left(M\right)$, the following properties
are equivalent : (i) $f$ is basic ; (ii) $f$ is constant on each
leaf of $\mathfrak{F}$.

We recall also, that a vector field $X\in\mathfrak{X}(M)$ is said
to be foliate, or that it is an infinitesimal automorphism of $\mathfrak{F}$
if in a neighborhood of any point of $M$, the local one parameter
group associated to $X$ leaves the foliation $\mathfrak{F}$ invariant. 

We have the following equivalence : (i) $X$ is foliate. (ii) $[X,Y]\in\Gamma(E)$
for all $Y\in\Gamma(E)$. (iii) In a local coordinate system $(x^{pi},y^{i})_{1\leq p\leq k,1\leq i\leq n}$,
the vector field $X$ has the following form 
\[
X=\sum_{p,i}\xi^{pi}\left(\left(x^{qj}\right)_{q,j},y^{1},\ldots,y^{n})\right)\frac{\partial}{\partial x^{pi}}+\sum_{j=1}^{n}\eta^{j}(y^{1},\ldots,y^{n})\frac{\partial}{\partial y_{j}}.
\]
We denote by$\mathfrak{\mathscr{I}}(M,\mathfrak{F})$ the space of
foliate vector fields for $\mathfrak{F}.$ The following properties
are satisfied :
\begin{enumerate}
\item $\mathfrak{\mathscr{I}}(M,\mathfrak{F})$ is a Lie algebra.
\item $\mathfrak{\mathscr{I}}(M,\mathfrak{F})$ is a module over the ring
$\mathscr{A}_{b}^{0}\left(M,\mathfrak{F}\right)$ of basic functions. 
\end{enumerate}
A smooth $r$-form $\alpha$ a on $M$ is said to be basic (\cite{key-MLN2})
if: $i(Y)\alpha=0\,\,\,\mathrm{et}\,\,\,i(Y)d\alpha=0\,\,\,\,\,\forall Y\in\Gamma(E).$ 

The following properties are equivalent : (i) $\alpha$ is basic.
(ii) In every simple distinguished open set, equipped with Darboux's
local coordinate system, $\alpha$ take the form 
\[
\alpha=\sum_{1\leq i_{1}<\ldots<i_{r}\leq n}\alpha_{i_{1}\ldots i_{r}}dy_{i_{1}}\land\ldots\land dy_{i_{r}}
\]
where the coefficients $\alpha_{i_{1}\ldots i_{r}}$ are basic functions.

\section{Polarized Hamiltonian vector fields}

The notations being the same as in the previous paragraph. 

Let $M$ be a manifold endowed with a polarized $k$-symplectic structure
$\left(\theta,\mathfrak{\mathit{E}}\right).$ Consider the linear
mapping $\zeta:\mathfrak{X}(M)\longrightarrow\mathscr{A}_{1}(M)\otimes\mathbb{R}^{k}$,
defined by: $\zeta\left(X\right)=i(X)\theta,$ for all $X\in\mathfrak{X}(M)$. 
\begin{defn}
A vector filed $X\in\mathfrak{X}(M)$ is said to be locally Hamiltonian
polarized if it satisfies the following conditions: (i) $X$ is foliate;
(ii) the $\mathbb{R}^{k}$-valued $1-$form $\zeta(X)$ is closed.
\end{defn}

We denoted by $H^{0}\left(M,\mathfrak{F}\right)$ the real vector
space of locally Hamiltonians polarized vector fields 
\[
H^{0}\left(M,\mathfrak{F}\right)=\left\{ X\in\mathfrak{\mathscr{I}}(M,\mathfrak{F})\mid d(\zeta\left(X\right))=0\right\} 
\]
 Let $X$ be a locally Hamiltonian polarized vector field. Thus, locally,
around each point $x\in M$, there is an open neighborhood $U$ of
$x$ in $M$, and an $\mathbb{R}^{k}-$smooth mapping $H\in\mathscr{C}^{\infty}\left(U\right)\otimes\mathbb{R}^{k}$
such that $\zeta\left(X\right)=-dH$.

With respect to a Darboux's local coordinate system $(x^{pi},y^{i})_{1\leq p\leq k,1\leq i\leq n}$,
defined on an open neighborhood $U$ of $M$, the equations of motion
of $X$ are given by

\[
\begin{cases}
\frac{dx^{pi}}{dt}=-\frac{\partial H^{p}}{\partial y^{i}}\\
\delta_{q}^{p}\frac{dy^{i}}{dt}=\frac{\partial H^{p}}{\partial x^{qi}}\\
\frac{\partial H^{p}}{\partial x^{qi}}\in\mathscr{A}_{b}^{0}(M,\mathfrak{F}).
\end{cases}
\]

(Hamilton's equation of locally polarized Hamiltonian vector field),
and, the mapping $H$ and $X$ have the following forms respectively
\[
H=\sum_{p=1}^{k}\left(\sum_{j=1}^{n}a_{j}(y^{1},...,y^{n})x^{pj}+b^{p}(y^{1},...,y^{n})\right)\otimes v_{p}
\]
 and 
\[
X_{H}=-\sum_{s=1}^{n}\sum_{p=1}^{k}\left(\sum_{j=1}^{n}x^{ps}\frac{\partial a_{j}}{\partial y^{s}}+\frac{\partial b^{p}}{\partial y^{s}}\right)\frac{\partial}{\partial x^{ps}}+\sum_{j=1}^{n}a_{j}\frac{\partial}{\partial y^{j}}
\]

where $a_{j},b^{p}\in\mathscr{A}_{b}^{0}\left(U,\mathfrak{F_{\mathit{U}}}\right).$
\begin{defn}
An element $X\in\mathfrak{\mathscr{I}}(M,\mathfrak{F})$ is said to
be a polarized Hamiltonian vector field if the $\mathbb{R}^{k}$-valued
one form $\zeta\left(X\right)$ is\textbf{ }exact. We denote by $H\left(M,\mathfrak{F}\right)$
the real vector space of polarized Hamiltonian vector fields.
\end{defn}

We assume that $M$ is connected. Then we have an exact sequence of
vector spaces 
\[
0\longrightarrow\mathbb{R}^{k}\overset{i}{\longrightarrow}\mathscr{C}^{\infty}\left(M,\mathbb{R}^{k}\right)\overset{d}{\longrightarrow}\mathscr{A}_{1}(M)\otimes\mathbb{R}^{k}\longrightarrow0.
\]

It is clear that $\zeta(H(M,\mathfrak{F}))$ is vector subspace of
$\mathscr{A}_{1}(M)\otimes\mathbb{R}^{k}$. Then we take 
\[
\mathfrak{H}(M,\mathfrak{F})=d^{-1}\left(\zeta(H(M,\mathfrak{F}))\right).
\]

Therefore, for every mapping $H\in\mathscr{C}^{\infty}\left(M\right)\otimes\mathbb{R}^{k}$
we have the following equivalence :
\begin{enumerate}
\item $H\in\mathfrak{H}(M,\mathfrak{F})$ ;
\item there exists a unique polarized Hamiltonian vector field $X_{H}\in H(M,\mathfrak{F})$
such that $\zeta(X_{H})=-dH.$ 
\end{enumerate}
The elements of $\mathfrak{H}(M,\mathfrak{F})$ are said to be polarized
Hamiltonians, and $X_{H}$, the corresponding polarized Hamiltonian
vector fields.

Also, we have the mapping $\nu:\mathfrak{H}(M,\mathfrak{F})\longrightarrow H(M,\mathfrak{F})$
defined by: $\nu\left(H\right)=X_{H}.$ 
\[
\]
\begin{prop}
The following diagram:
\[
\begin{array}{ccc}
H(M,\mathfrak{F}) & \overset{\zeta}{\longrightarrow} & \mathscr{A}^{1}(M)\otimes\mathbb{R}^{k}\\
\nwarrow &  & \nearrow\\
\nu & \mathfrak{H}(M,\mathfrak{F}) & -d
\end{array}
\]
is commutative 
\end{prop}

\section{Polarized vectorial Poisson structure subordinate to polarized $k$-symplectic
structure }

The hypothesis and notations being the same as above.

Recall that (\cite{key-AW}, \cite{key-AW-GZ}), if $H,K\mathfrak{\in H}(M,\mathfrak{F})$
with associated polarized Hamiltonian vector fields $X_{H}$ , $X_{K}$
respectively, then the Lie bracket $[X_{H},X_{K}]$ is a polarized
Hamiltonian vector field which is associated to polarized Hamiltonian
$\left\{ H,K\right\} $ defined by 
\[
\left\{ H,K\right\} =\left\{ H,K\right\} ^{p}\otimes v_{p}=-\theta^{p}\left(X_{H},X_{K}\right)\otimes v_{p}
\]
 i.e. $[X_{H},X_{K}]=X_{\left\{ H,K\right\} }$and the correspondence
$(H,K)\longmapsto\left\{ H,K\right\} $ from $\mathfrak{H}(M,\mathfrak{F})\times\mathfrak{H}(M,\mathfrak{F})$
into $\mathfrak{H}(M,\mathfrak{F})$, gives to $\mathfrak{H}(M,\mathfrak{F})$
a structure of Lie algebra.
\begin{defn}
The Lie algebra $\left(\mathfrak{H}(M,\mathfrak{F}),\left\{ ,\right\} \right)$
is said to be the polarized Poisson structure subordinate to the polarized
$k$-symplectic structure $\left(\theta,E\right)$.
\end{defn}

\begin{prop}
We have the following properties:
\begin{enumerate}
\item $H\left(M,\mathfrak{F}\right)$ is a real Lie algebra.
\item $\left[H^{0}\left(M,\mathfrak{F}\right),H^{0}\left(M,\mathfrak{F}\right)\right]\subset H\left(M,\mathfrak{F}\right)$. 
\item $H\left(M,\mathfrak{F}\right)$ is an ideal of $H^{0}\left(M,\mathfrak{F}\right)$.
\item The sequence of Lie algebras 
$0\longrightarrow\mathbb{R}^{k}\longrightarrow\mathfrak{H}(M,\mathfrak{F})\longrightarrow H\left(M,\mathfrak{F}\right)\hookrightarrow H^{0}(M,\mathfrak{F})\longrightarrow\frac{H^{0}\left(M,\mathfrak{F}\right)}{H(M,\mathfrak{F})}\longrightarrow0$ is exact. 
\end{enumerate}
\end{prop}

Locally, with respect to a Darboux's local coordinate system $(x^{pi},y^{i})_{1\leq p\leq k,1\leq i\leq n}$,
defined on an open neighborhood $U$ of $M$, the bracket $\left\{ H,K\right\} $
is written in the form $\left\{ H,K\right\} ^{p}  \otimes v_{p}$ where 

\begin{align*}
\left\{ H,K\right\} ^{p} &=\sum_{i=1}^{n}\left(\frac{\partial H^{p}}{\partial y^{i}}\frac{\partial K^{p}}{\partial x^{pi}}-\frac{\partial H^{p}}{\partial x^{pi}}\frac{\partial K^{p}}{\partial y^{i}}\right)\\
\\
 & =\left(\frac{\partial}{\partial y^{i}}\otimes v^{p}\right)\wedge\left(\frac{\partial}{\partial x^{pi}}\otimes v^{p}\right)\left(dH^{l}\otimes v_{l},dK^{r}\otimes v_{r}\right)
\end{align*}
where $\left(v^{p}\right)_{1\leq p\leq k}$ is the dual basis of the
standard basis $\left(v_{p}\right)_{1\leq p\leq k}$ of $\mathbb{R^{\mathit{k}}}$.

By using the Poisson polarized bracket, the Hamilton's equations of
polarized Hamiltonian vector field $X_{H}$ take the new form

\[
\begin{cases}
\frac{dx^{pi}}{dt}=\left\{ \left(x^{1i},\ldots,x^{ki}\right),H\right\} ^{p}\\
\delta_{q}^{p}\frac{dy^{i}}{dt}=\left\{ \left(y^{i}\delta_{1}^{p},\ldots,y^{i}\delta_{q}^{p},\ldots,y^{i}\delta_{k}^{p}\right),H\right\} ^{p}\\
\frac{\partial H^{p}}{\partial x^{qi}}\in\mathscr{A}_{b}^{0}(M,\mathfrak{F}).
\end{cases}
\]

Therefore, for every point $x$ of $U$ we have 
\[
\begin{array}{ccccc}
\left\{ H,K\right\} (x) & = & \left(X_{K}\cdot H\right)(x) & = & \left\langle dH(x),X_{K}(x)\right\rangle \\
 & = & -\left(X_{H}\cdot K\right)(x) & =- & \left\langle dK(x),X_{H}(x)\right\rangle 
\end{array}
\]
And, Let us assume that the function $H$ is fixed in $\mathfrak{H}\left(M,\mathfrak{F}\right)$
and when $K$ varies, $\left\{ H,K\right\} (x)$ only depends upon
$dK(x)$. Similarly, when we assume that $K$ is fixed in $\mathfrak{H}\left(M,\mathfrak{F}\right)$,
we can show that when $H$ varies, $\left\{ H,K\right\} (x)$ only
depends upon $dH(x)$. Therefore, for every point $x$ of $M$ there
exists a bilinear, skew symmetric mapping 
\[
P(x):\left(T_{x}^{*}M\otimes\mathbb{R^{\mathit{k}}}\right)\times\left(T_{x}^{*}M\otimes\mathbb{R^{\mathit{k}}}\right)\longrightarrow\mathbb{R^{\mathit{k}}}
\]
such that $P(x)\left(dH(x),dK(x)\right)=\left\{ H,K\right\} (x)$
for all $H$ and $K$ in $\mathfrak{H}\left(M,\mathfrak{F}\right)$,
$P$ is called polarized vectorial Poisson tensor.\textbf{ }

With respect to the local coordinate system, $(x^{pi},y^{i})_{1\leq p\leq k,1\leq i\leq n}$,
defined on an open neighborhood $U$ of $M$, the vectorial tensor
$P$ is

\[
P=\sum_{p=1}^{k}\sum_{i=1}^{n}\left(\left(\frac{\partial}{\partial y^{i}}\otimes v^{l}\right)\wedge\left(\frac{\partial}{\partial x^{pi}}\otimes v^{l}\right)\right)\otimes v_{p}
\]
\begin{rem}
The vectorial tensor $P$ vanishes on the annihilator $\mathscr{A}_{b}^{1}(M)\otimes\mathtt{\mathbb{R}}^{\mathit{k}}$
of $E$ in the space $\mathscr{A}^{1}(M)\otimes\mathtt{\mathbb{R}}^{\mathit{k}}.$
\end{rem}

\section{Polarized vectorial Poisson manifolds}

Let $\left(M,\mathfrak{F}\right)$ be a foliated manifold.
\begin{defn}
A polarized Poisson structure on $\left(M,\mathfrak{F}\right)$ is
a pair $\left(\mathfrak{H}\left(M,\mathfrak{F}\right),\left\{ ,\right\} \right)$
in which $\mathfrak{H}\left(M,\mathfrak{F}\right)$ is a vector subspace
of $\mathscr{C}^{\infty}\left(M\right)\otimes\mathtt{\mathbb{R}}^{\mathit{k}}.$
and $\left\{ ,\right\} $ is an $\mathbb{R\mathit{-}}$ bilinear mapping
from $\mathfrak{H}\left(M,\mathfrak{F}\right)\times\mathfrak{H}\left(M,\mathfrak{F}\right)$into
$\mathfrak{H}\left(M,\mathfrak{F}\right)$ satisfying: 
\begin{enumerate}
\item $\left(\mathfrak{H}\left(M,\mathfrak{F}\right),\left\{ ,\right\} \right)$
is a Lie algebra. 
\item $\left\{ H,K\right\} =0,\,\,\,\forall H,K\in\mathscr{A}_{b}^{0}\left(M\right)\otimes\mathtt{\mathbb{R}}^{\mathit{k}}$. 
\item For each $H\in\mathfrak{H}\left(M,\mathfrak{F}\right)$, there is
a foliate vector field $X_{H}$ on $M$ such that 
\[
X_{H}\left(K\right)=\left\{ K,H\right\} ,\,\forall K\in\mathfrak{H}\left(M,\mathfrak{F}\right).
\]
\end{enumerate}
\end{defn}
A polarized vectorial Poisson structure on $\left(M,\mathfrak{F}\right)$
can be defined on $M$ by a $\mathscr{C}^{\infty}(M)$\textminus bilinear
skew symmetric map
\[
P:\mathscr{A}^{1}(M)\otimes\mathtt{\mathbb{R}}^{\mathit{k}}\times\mathscr{A}^{1}(M)\otimes\mathtt{\mathbb{R}}^{\mathit{k}}\longrightarrow\mathscr{C}^{\infty}\left(M,\mathtt{\mathbb{R}}^{\mathit{k}}\right)
\]
satisfying the following properties: 
\begin{enumerate}
\item $\forall H,K\in\mathfrak{H}(M,\mathfrak{F})$, $P(dH,dK)\in\mathfrak{H}(M,\mathfrak{F})$.
\item The correspondence $(H,K)\longmapsto P(dH,dK)=\left\{ H,K\right\} $
from $\mathfrak{H}(M,\mathfrak{F})\times\mathfrak{H}(M,\mathfrak{F})$
into $\mathfrak{H}(M,\mathfrak{F})$, gives to $\mathfrak{H}(M,\mathfrak{F})$
a structure of Lie algebra.
\item $P$ vanishes on $\mathscr{A}_{b}^{1}(M)\otimes\mathtt{\mathbb{R}}^{\mathit{k}}$:
annihilator of $E$ in the space $\mathscr{A}^{1}(M)\otimes\mathtt{\mathbb{R}}^{\mathit{k}}.$
\item For every $H\in\mathfrak{H}\left(M,\mathfrak{F}\right)$ there is
a foliate vector field $X_{H}$ such that 
\[
P(dH,dK)=-X_{H}(K),\,\,\,\forall K\in\mathfrak{H}(M,\mathfrak{F}).
\]
\end{enumerate}

\section{Some properties of the polarized vectorial Poisson structures}

Let $\left(\mathfrak{H}\left(M,\mathfrak{F}\right),P\right)$ be a
polarized vectorial Poisson structure on a foliated manifold $\left(M,\mathfrak{F}\right)$.
Let $(x_{1},...,x_{p},y_{1},\ldots,y_{q})$ be distinguished coordinates
of a local foliate chart. Locally the polarized vectorial tensor $P$
take the form 

\[
\begin{array}{ccc}
P & = & A_{ls}^{ijr}\left(\left(\frac{\partial}{\partial x_{i}}\otimes v^{l}\right)\wedge\left(\frac{\partial}{\partial y_{j}}\otimes v^{s}\right)\right)\otimes v_{r}\\
 & + & B_{ls}^{ijr}\left(\left(\frac{\partial}{\partial x_{i}}\otimes v^{l}\right)\wedge\left(\frac{\partial}{\partial x_{j}}\otimes v^{s}\right)\right)\otimes v_{r}
\end{array}
\]
And for all $\alpha\in\mathscr{A}_{1}(M)\otimes\mathtt{\mathbb{R}}^{\mathit{k}},$
we associate a $\mathscr{C}^{\infty}(M)-$linear mapping
\[
P(\alpha,\cdot):\mathscr{A}_{1}(M)\otimes\mathtt{\mathbb{R}}^{\mathit{k}}\longrightarrow\mathscr{C}^{\infty}\left(M\right)\otimes\mathtt{\mathbb{R}}^{\mathit{k}}
\]
 such that $P(\alpha,\cdot)(\beta)=P(\alpha,\beta)$, and a mapping
\[
\Xi:\mathfrak{X}(M)\longrightarrow\mathcal{L}_{\mathscr{C}^{\infty}\left(M\right)}\left(\mathscr{A}_{1}(M)\otimes\mathtt{\mathbb{R}}^{\mathit{k}},\mathscr{C}^{\infty}\left(M\right)\otimes\mathtt{\mathbb{R}}^{\mathit{k}}\right)
\]

that we define by
\[
\Xi\left(X\right)\left(\beta\right)=\left\langle \beta,X\right\rangle =\sum_{p=1}^{k}\beta^{p}(X)v_{p}=\sum_{p=1}^{k}\left(\beta^{p}\otimes v_{p}\right)(X).
\]

Locally, with respect to distinguished coordinates $(x_{1},...,x_{p},y_{1},\ldots,y_{q})$
we have :
\[
\Xi\left(\frac{\partial}{\partial x^{l}}\right)\left(\beta\right)=\sum_{p=1}^{k}\frac{\partial\beta^{p}}{\partial x^{l}}v_{p}
\]
\begin{itemize}
\item The mapping $\Xi$ is an isomorphism if and only if $k=1$. 
\item In particular $\Xi\left(X_{H}\right)\left(dK\right)=-P(dH,dK)=-\left\{ H,K\right\} =\left\{ K,H\right\} $
for all $H,\,K\in\mathfrak{H}(M,\mathfrak{F})$. 
\end{itemize}

\section{Polarized vectorial Poisson structure subordinate to the natural
polarization $k-$symplectic structure of $\hom\left(\mathbb{\mathcal{G}\mathit{,}R^{\mathrm{\mathit{k+1}}}}\right)$ }

Let $\mathcal{\left(\mathcal{G}\mathit{,\left[,\right]}\right)}$
be a real $n$-dimensional Lie algebra equipped with a basis $\left(e_{i}\right)_{1\leq i\leq n}$.
We denote by $\left(\omega^{i}\right)_{1\leq i\leq n}$ its dual basis
and $C_{ij}^{k}$ the structure constants of $\mathcal{G}$: $\left[e_{i},e_{j}\right]=C_{ij}^{k}e_{k}.$
And let
\[
\hom\left(\mathbb{\mathcal{G}\mathit{,}R^{\mathit{k+1}}}\right)=\mathbb{\mathcal{G}^{\mathit{*}}\otimes R^{\mathit{k+1}}}
\]
the space of linear maps of $\mathcal{G}$ into $\mathbb{R}^{\mathit{k+1}}$.
The vector space $\hom\left(\mathbb{\mathcal{G}\mathit{,}R^{\mathit{k+1}}}\right)$
is generated by the linear maps:
\[
\omega^{i}\otimes v_{q},\,\omega^{i}\otimes w\,\,\,\left(1\leq q\leq k;1\leq i\leq n\right)
\]
where $\left(v_{q},w\right)_{1\leq q\leq k;}$ is the canonical basis
of $\mathbb{R}^{k+1}$.

Each element $X$ of $\hom\left(\mathbb{\mathcal{G}\mathit{,}\mathbb{R}^{\mathrm{\mathit{k+1}}}}\right)$
can be written as
\[
X=x_{i}^{q}\omega^{i}\otimes v_{q}+y_{i}\omega^{i}\otimes w=\left(\begin{array}{ccc}
x_{1}^{1} & \ldots & x_{n}^{1}\\
\vdots & \vdots & \vdots\\
x_{1}^{k} & \vdots & x_{n}^{k}\\
y_{1} & \cdots & y_{n}
\end{array}\right).
\]
The linear mapping $X:\mathcal{G}\mathit{\longrightarrow}\mathbb{R}^{\mathrm{\mathit{k+1}}}$
transforms $u=u^{j}e_{j}$ into
\[
\begin{array}{ccc}
X\left(u\right) & = & x_{i}^{q}\omega^{i}\otimes v_{q}\left(u^{j}e_{j}\right)+y_{i}\omega^{i}\otimes w\left(u^{j}e_{j}\right)\\
 & = & x_{i}^{q}u^{i}v_{q}+y_{i}u^{i}w
\end{array}
\]
So, in terms of matrices, we have
\[
X\left(u\right)=\left(\begin{array}{ccc}
x_{1}^{1} & \ldots & x_{n}^{1}\\
\vdots & \vdots & \vdots\\
x_{1}^{k} & \vdots & x_{n}^{k}\\
y_{1} & \cdots & y_{n}
\end{array}\right)\left(\begin{array}{c}
u^{1}\\
\vdots\\
u^{n}
\end{array}\right).
\]
 The space $\hom\left(\mathbb{\mathcal{G}\mathit{,}R^{\mathit{\mathrm{\mathit{k+1}}}}}\right)$
is an $n\left(k+1\right)$-dimensional smooth manifold. We equip this
space with the coordinate system $\left(x_{i}^{q},y_{i}\right)_{1\leq q\leq k;1\leq i\leq n}$.

There is a natural polarized $k$-symplectic structure $\left(\theta,\mathfrak{F}\right)$
on $\hom\left(\mathbb{\mathcal{G}\mathit{,}R^{\mathit{\mathrm{\mathit{k+1}}}}}\right)$
defined by 
\[
\theta=\left(\sum_{i=1}^{n}dx_{i}^{q}\wedge dy_{i}\right)\otimes v_{q}
\]
 and $\mathfrak{F}$ is the foliation defined by $dy_{1}=0,\cdots,dy_{n}=0$.

The associated polarized Hamiltonians are the $\mathbb{R^{\mathit{k}}}-$valued
smooth functions 
\[
H\in\mathscr{C}^{\infty}\left(\hom\left(\mathbb{\mathcal{G}\mathit{,}R^{\mathit{k+1}}}\right)\right)\otimes\mathbb{R}^{\mathit{k}}
\]
 defined on $\hom\left(\mathbb{\mathcal{G}\mathit{,}R^{\mathit{k+1}}}\right)$
by the following expressions
\[
H\left(X\right)=\left(a^{i}\left(y_{1},\cdots,y_{n}\right)x_{i}^{q}+b^{q}\left(y_{1},\cdots,y_{n}\right)\right)\otimes v_{q}
\]
where $a^{1},\ldots,a^{n},b^{q}\left(q=1,\cdots,k\right)$ are real
basic functions. 

The polarized Poisson bracket of $H\left(X\right)=\left(a^{i}\left(y_{1},\cdots,y_{n}\right)x_{i}^{q}+b^{q}\left(y_{1},\cdots,y_{n}\right)\right)\otimes v_{q}$
and $K\left(X\right)=\left(a'^{i}\left(y_{1},\cdots,y_{n}\right)x_{i}^{q}+b'^{q}\left(y_{1},\cdots,y_{n}\right)\right)\otimes v_{q}$
 is given by
\[
\begin{array}{ccc}
\left\{ H,K\right\}^{q} \left(X\right) & = & \left(\sum_{i=1}^{n}\left(\frac{\partial H^{q}}{\partial y^{i}}\frac{\partial K^{q}}{\partial x_{i}^{q}}-\frac{\partial H^{q}}{\partial x_{i}^{q}}\frac{\partial K^{q}}{\partial y^{i}}\right)\right)\\
 & = & \left(\sum_{i=1}^{n}\left(\left(x_{j}^{q}\frac{\partial a^{j}}{\partial y^{i}}+\frac{\partial b^{q}}{\partial y^{i}}\right)a'^{i}-a_{i}\left(x_{j}^{q}\frac{\partial a'^{j}}{\partial y^{i}}+\frac{\partial b'^{q}}{\partial y^{i}}\right)\right)\right)\\
 & = & \left(\sum_{i=1}^{n}\left(x_{j}^{q}\left(a'^{i}\frac{\partial a^{j}}{\partial y^{i}}-a^{i}\frac{\partial a'^{j}}{\partial y^{i}}\right)+\left(a'^{i}\frac{\partial b^{q}}{\partial y^{i}}-a^{i}\frac{\partial b'^{q}}{\partial y^{i}}\right)\right)\right).
\end{array}
\]

The bracket, so defined, allows to provide $\mathcal{\mathfrak{H}}\left(\hom\left(\mathbb{\mathcal{G}\mathit{,}R^{\mathit{k+1}}}\right),\mathfrak{F}\right)$
with a polarized vectorial Poisson structure subordinate to the polarization
$k$-symplectic structure $\left(\theta,\mathfrak{F}\right).$ This
structure does not depend on the law of $\mathcal{G}$.

\section{Linear polarized vectorial Poisson structure of $\hom\left(\mathbb{\mathcal{G}\mathit{,}R^{\mathit{k+1}}}\right)$}

In addition to the polarized vectorial Poisson structure subordinate
to the natural $k-$symplectic polarization on $\hom\left(\mathbb{\mathcal{G}\mathit{,}R^{\mathit{k+1}}}\right),$
we can define another polarized vectorial Poisson structure $\left(\mathcal{\mathfrak{H}}\left(\hom\left(\mathbb{\mathcal{G}\mathit{,}R^{\mathit{k+1}}}\right),\mathfrak{F}\right);\left\{ ,\right\} ^{L}\right)$
so-called the linear polarized vector Poisson structure of $\left(\mathcal{G},\left[,\right]\right).$ 

Let $H\in\mathcal{\mathfrak{H}}\left(\hom\left(\mathbb{\mathcal{G}\mathit{,}R^{\mathit{k+1}}}\right),\mathfrak{F}\right)$,
$X\in\hom\left(\mathbb{\mathcal{G}\mathit{,}R^{\mathit{k+1}}}\right)$
and $j_{q}:\mathcal{G}^{*}\longrightarrow\hom\left(\mathbb{\mathcal{G}\mathit{,}R^{\mathit{k+1}}}\right)$
defined by
\[
j_{q}(\omega^{i})=\omega^{i}\otimes v_{q}.
\]
 For all $p,q\in\left\{ 1,\cdots,k\right\} ,$ the composite map
\[
\mathcal{G^{\mathrm{*}}}\stackrel{j_{q}}{\longrightarrow}\hom\left(\mathbb{\mathcal{G}\mathit{,}R^{\mathit{k+1}}}\right)\stackrel{dH_{X}^{p}}{\longrightarrow}\mathbb{R^{\mathit{}}}
\]
is the linear form on $\mathcal{G^{\mathit{*}}}$ defined by 
\[
\left(dH_{X}^{p}\circ j_{q}\right)\left(\omega^{i}\right)=dH_{X}^{p}\left(\omega^{i}\otimes v_{q}\right)=\frac{\partial H^{p}}{\partial x_{i}^{q}}\left(X\right)=\delta_{q}^{p}a^{i}.
\]
Hence, 
\[
\left.\left(dH_{X}^{1}\circ j_{1}\right)=\cdots=\left(dH_{X}^{k}\circ j_{k}\right)=a^{i}e_{i}.\right.
\]
We take 
\[
\begin{array}{ccc}
\left\{ H,K\right\} ^{L}\left(X\right) & = & pr_{\mathbb{R^{\mathit{k}}}}\left\langle \left[dH_{X}^{1}\circ j_{1},dK_{X}^{1}\circ j_{1}\right],X\right\rangle \\
 & = & \cdots\\
 & = & pr_{\mathbb{R^{\mathit{k}}}}\left\langle \left[dH_{X}^{k}\circ j_{k},dK_{X}^{k}\circ j_{k}\right],X\right\rangle 
\end{array}
\]
 therefore, 
\[
\begin{array}{ccc}
\left\{ H,K\right\} ^{L}\left(X\right) & = & pr_{\mathbb{R^{\mathit{k}}}}\left\langle \left[dH_{X}^{p}\circ j_{p},dK_{X}^{p}\circ j_{p}\right],X\right\rangle \\
 & = & pr_{\mathbb{R^{\mathit{k}}}}\left\langle \left[a^{i}e_{i},a'^{j}e_{j}\right],X\right\rangle \\
 & = & pr_{\mathbb{R^{\mathit{k}}}}\left\langle \left(a^{i}a'^{j}C_{ij}^{l}e_{l}\right),X\right\rangle \\
 & = & \left(a^{i}a'^{j}C_{ij}^{l}x_{l}^{p}\right)v_{p}\\
 & = & \left(\sum_{1\leq i<j\leq n}C_{ij}^{l}\left(a^{i}a'^{j}-a^{j}a'^{i}\right)x_{l}^{p}\right)v_{p}
\end{array}
\]

\subsection*{Examples}
\begin{enumerate}
\item $\mathcal{G}$ is an abelian Lie algebra. In this case $\left\{ ,\right\} ^{L}=0$.
Consequently, $\left(\mathcal{\mathfrak{H}}\left(\hom\left(\mathbb{\mathcal{G}\mathit{,}R^{\mathit{k+1}}}\right),\mathfrak{F}\right);\left\{ ,\right\} ^{L}\right)$
the abelian polarized vectorial Poisson.
\item $\mathcal{G}$ is Heisenberg\textquoteright s Lie algebra $\mathcal{H}_{1}$
of dimension 3. The Lie algebra law of $\mathcal{H}_{1}$ is given
by $[e_{1},e_{2}]=e_{3}$. And so, for all $H,K\in\mathcal{\mathfrak{H}}\left(\hom\left(\mathbb{\mathcal{G}\mathit{,}R}^{k+1}\right),\mathfrak{F}\right)$,
$X\in\hom\left(\mathbb{\mathcal{G}\mathit{,}R^{\mathit{k+1}}}\right)$,
where $H(X)=\left(a^{i}\left(y_{1},y_{2},y_{3}\right)x_{i}+b\left(y_{1},y_{2},y_{3}\right)\right)\otimes v_{q}$
and $K(X)=\left(a'^{i}\left(y_{1},y_{2},y_{3}\right)x_{i}+b'\left(y_{1},y_{2},y_{3}\right)\right)\otimes v_{q}$,
we have 
\[
\left\{ H,K\right\} ^{L}(X)=\left(\left(a^{1}a'^{2}-a^{2}a'^{1}\right)x_{3}^{p}\right)v_{p}.
\]
\item $\mathcal{G}$ is the polarized 1-symplectic nilpotent Lie algebra
$h_{3}\oplus a$ of dimension 4. The Lie algebra $\mathcal{G}$ defined
by
\begin{align*}
d\omega_{1} & =\omega_{2}\wedge\omega_{3}\\
d\omega_{i} & =0\,\,\,\,\,\,(i=2,3,4).
\end{align*}
So, for all $H,K\in\mathcal{\mathfrak{H}}\left(\hom\left(\mathbb{\mathcal{G}\mathit{,}R^{\mathit{k+1}}}\right),\mathfrak{F}\right)$,
$X\in\hom\left(\mathbb{\mathcal{G}\mathit{,}R^{\mathit{k+1}}}\right)$,
where $H(X)=\left(a^{i}\left(y_{1},y_{2},y_{3},y_{4}\right)x_{i}^{q}+b^{q}\left(y_{1},y_{2},y_{3},y_{4}\right)\right)\otimes v_{q}$and

$K(X)=\left(a'^{i}\left(y_{1},y_{2},y_{3},y_{4}\right)x_{i}^{q}+b'^{q}\left(y_{1},y_{2},y_{3},y_{4}\right)\right)\otimes v_{q}$,
we have 
\[
\left\{ H,K\right\} ^{L}(X)=\left(\left(a^{3}a'^{2}-a^{2}a'^{3}\right)x_{1}^{p}\right)v_{p}.
\]
\item $\mathcal{G}$ is the polarized 1-symplectic nilpotent Lie algebra
$n_{4}$ of dimension 4. The Lie algebra $\mathcal{G}$ defined by
\begin{align*}
d\omega_{3} & =\omega_{1}\wedge\omega_{2}\\
d\omega_{4} & =\omega_{1}\wedge\omega_{3}\\
d\omega_{i} & =0\,\,\,\,\,\,(i=1,2).
\end{align*}
So, for all $H,K\in\mathcal{\mathfrak{H}}\left(\hom\left(\mathbb{\mathcal{G}\mathit{,}R^{\mathit{k+1}}}\right),\mathfrak{F}\right)$,
$X\in\hom\left(\mathbb{\mathcal{G}\mathit{,}R^{\mathit{k+1}}}\right)$,
where $H(X)=\left(a^{i}\left(y_{1},y_{2},y_{3},y_{4}\right)x_{i}^{q}+b^{q}\left(y_{1},y_{2},y_{3},y_{4}\right)\right)\otimes v_{q}$
and

$K(X)=\left(a'^{i}\left(y_{1},y_{2},y_{3},y_{4}\right)x_{i}^{q}+b'^{q}\left(y_{1},y_{2},y_{3},y_{4}\right)\right)\otimes v_{q}$,
we have 
\[
\left\{ H,K\right\} ^{L}(X)=\left(\left(a^{2}a'^{1}-a^{1}a'^{2}\right)x_{3}^{p}\right)v_{p}+\left(\left(a^{3}a'^{1}-a^{1}a'^{3}\right)x_{4}^{p}\right)v_{p}.
\]
\end{enumerate}

\end{document}